\documentclass{amsart}

\usepackage{amsmath,amssymb,amsthm,a4wide}
\usepackage{graphicx,tikz}

\newtheorem*{thm}{Theorem}
\newtheorem*{proposition}{Proposition}

\theoremstyle{definition}

\theoremstyle{remark}

\begin{document}

\title[]{ Numerical Integration on Graphs: \\ where to sample and how to weigh}
\keywords{Graph, Sampling, Graph Laplacian, Sampling, Heat Kernel, Packing.}
\subjclass[2010]{05C50, 05C70, 35P05, 65D32}

\author[]{George C. Linderman}
\address{Program in Applied Mathematics, Yale University, New Haven, CT 06511, USA}
\email{george.linderman@yale.edu}

\author[]{Stefan Steinerberger}
\address{Department of Mathematics, Yale University, New Haven, CT 06511, USA}
\email{stefan.steinerberger@yale.edu}

\begin{abstract} Let $G=(V,E,w)$ be a finite, connected graph with weighted edges. We are interested in the problem of finding
a subset $W \subset V$ of vertices and weights $a_w$ such that
$$    \frac{1}{|V|}\sum_{v \in V}^{}{f(v)}  \sim \sum_{w \in W}{a_w f(w)}$$
for functions $f:V \rightarrow \mathbb{R}$ that are `smooth' with respect to the geometry of the graph. The main
application are problems where $f$ is known to somehow depend on the underlying graph but is expensive to evaluate
on even a single vertex. We prove an inequality showing that the integration problem can be rewritten
as a geometric problem (`the optimal packing of heat balls'). We discuss how one would construct approximate solutions of the heat ball packing
problem; numerical examples demonstrate the efficiency of the method.

\end{abstract}

\maketitle

\section{Introduction}
\subsection{Introduction.} The purpose of this paper is to report on a very general idea in the context of sampling on graphs. It is a companion paper to \cite{stein1} dealing with a problem in Spectral Graph Theory. We extend some of the ideas from \cite{stein, stein1} to sampling, prove an inequality bounding the integration error in terms of the geometry of the sampling points and give several examples.
We will, throughout the paper, use $G=(V,E,w)$ to denote a connected graph with weighted
edges. 

\begin{quote}
\textbf{Problem} (Quadrature)\textbf{.} If we are allowed to sample in a set $W$ with $|W| =k$ vertices, which vertices and weights $a_w$ should we pick so that
$$    \frac{1}{|V|}\sum_{v \in V}^{}{f(v)}  \sim \sum_{w \in W}{a_w f(w)}$$
for functions $f:V \rightarrow \mathbb{R}$ that are `smooth with respect to to the geometry' of $G$?
\end{quote}
In many cases, this question does not make a lot of sense: since the graph is
finite, one can simply compute the true average of $f$ by summing over all $n$
vertices. This question is only interesting whenever sampling $f$ is difficult
or expensive. A toy example is the following: suppose we have a medical
database of $n \gg 1$ people containing all sorts of information and are
interested in the average blood pressure (not in the database) of those $n$ people. Actually going
out and measuring the blood pressure of all $n$ people would take a very long
time. However, it is known that blood pressure is strongly correlated with some
of the factors we have on file (say, age, weight, and smoking habits) and weakly
connected or unconnected to others (say, eye color). We can then build a
weighted graph on $n$ vertices where the weight on the edge connecting two people depends
on how similar they are with regards to relevant factors --  the hope is that
blood pressure, as a function on the graph, is then smoothly varying. Which of
the, say, $n/1000$ people should we do a blood pressure measurement on so that the
sample average is representative of the global average? Is there a way to pick
them in a way that decreases the expected error over a pure random selection?
It is likely possible to make stronger statements if one restricts to certain
classes of graphs and functions and this could be of substantial interest; our paper will
only address the most general case.

\subsection{Formal Setup.} We will now make these notions precise. As is not
surprising, the precise definition of `smoothness' of the function is crucial.
If $f$ has no particular structure, then there is little hope of being able to
achieve anything at all. However, a weighted graph does induce a natural notion
of smoothness: we want the function to vary little over edges with a large
weight (since vertices connected by a large weight should be `similar') whereas
a large variation over edges with a smaller weight would not be that surprising.
We now introduce a notion of a Laplacian on a Graph; it slightly deviates from more classical
notions and could be of independent interest (see Section \S 2.4 for a more detailed discussion). 
 Let $A$ denote the
(weighted) adjacency matrix of $G$
$$ A_{ij} = \begin{cases} w(e_{ij}) \qquad &\mbox{if}~i \sim_{E} j \\
0  &\mbox{otherwise,} \end{cases}
$$
where $w(e_{ij}) \geq 0$ is the weight of the edge $e_{ij} = e_{ji}$ connecting vertices $i$ and $j$.
$A$ is a symmetric matrix. We introduce the maximum sum of any of the rows of this symmetric matrix (coinciding, by symmetry,
with the maximum sum of any of the columns)
and use it to define a normalized adjacency matrix: more precisely, we have
$$d_{{\small \mbox{max}}} = \max_{1 \leq i \leq n} \sum_{j=1}^n A_{ij} \qquad \mbox{and set} \qquad A' = \frac{1}{d_{\small \mbox{max}}} A.$$
Finally, we introduce the (diagonal) degree matrix $D'$ associated to the renormalized adjacency matrix and use it to
define a Laplacian: we set
$$D'_{ii}  = \frac{1}{d_{\small \mbox{max}}}\sum_{j=1}^n A_{ij} \quad \mbox{and define} \quad L =  A' - D'.$$
We will never work directly
with the Laplacian: our main object of interest is the associated diffusion process whose generator
is given by
$$ P = L + \mbox{Id}_{n \times n},$$
where $\mbox{Id}_{n \times n}$ is the identity matrix of size $n$.
  $P$ is a symmetric
stochastic matrix and represents a lazy random walk where the probability of
``staying put'' depends on the vertex (as opposed to being, say,
$0.5$ as in the classical lazy random walk). We denote the eigenvalues of $P$, which are merely the eigenvalues of the Laplacian $L$ shifted by 1, by
 $\lambda_1, \dots, \lambda_n$. Since $P$ is a stochastic matrix, we have
$|\lambda_i| \leq 1$. The eigenvectors whose eigenvalues are close to 0 `diffuse quickly' and are thus the natural high-frequency
objects on the graph. This motivates an ordering of eigenvalues from low frequency to high frequency
$$1= |\lambda_1| \geq |\lambda_2| \geq |\lambda_3|\geq  \dots \geq |\lambda_n| \geq 0.$$

We denote the corresponding orthogonal eigenvectors of $P$, which clearly coincide with the eigenvectors of the Laplacian $L$, by $\phi_1, \dots, \phi_n \in
\mathbb{R}^n$ (normalized in $\ell^2$, $\phi_1
= 1/\sqrt{n}$ is the constant vector). We define a function space $X_{\lambda}$, the canonical
analogue of trigonometric polynomials on $\mathbb{T}^d$ or spherical harmonics on $\mathbb{S}^{d}$, via 
$$ X_{\lambda} = \left\{f: V \rightarrow \mathbb{R} \bigg| ~f = \sum_{|\lambda_{k}| \geq \lambda}{\left\langle f, \phi_k \right\rangle}\right\} \qquad \mbox{with norm} \qquad
 \|f\|^2_{X_{\lambda}} = \sum_{|\lambda_k| \geq \lambda}{ \left|\left\langle f, \phi_k \right\rangle\right|^2},$$
where $0 \leq
\lambda \leq 1$ is a parameter controlling the degree of smoothness. If
$\lambda > \mu$, then $X_{\lambda} \subset X_{\mu}$. Moreover, $X_{0} =
\left\{f:V \rightarrow \mathbb{R}\right\}$ contains all functions while, at
least on generic graphs, $X_1$ contains only the constant functions (depends on
whether the induced random walk is ergodic; this is not important for our
purposes). The norm is just the classical $L^2-$norm on the subspace -- in more
classical terms, we are simply considering the $L^2-$space obtained
via a Littlewood-Paley projection.  This function space $X_{\lambda}$ is
natural: if the graph $G$ approximates a torus $\mathbb{T}^d$, then this
function space will indeed approximate trigonometric functions. If $G$ is close
to a discretization of a sphere $\mathbb{S}^d$, then the space $X_{\lambda}$
approximates the space of low-degree spherical harmonics.

\section{Main Results}
\subsection{Main Result.} Our main result bounds the integration error in terms of $\|f\|_{X_{\lambda}}$ and a purely geometric quantity (explained in detail below) formulated in terms of the quadrature points and independent of $f$. This has the nice effect of multiplicatively separating the size of the function $\|f\|_{X_{\lambda}}$ and a quantity that can be interpreted as the 'quality' of the quadrature scheme. 
\begin{thm} 
\label{thm:powerits}
Let $W \subset V$ be equipped with weights $a_w$ summing to 1.  Then, for $\ell \in \mathbb{N}$, $0 < \lambda < 1$,
$$ \forall~f \in X_{\lambda} \quad   \left| \frac{1}{n} \sum_{v \in V}{f(v)} - \sum_{w \in W}{a_w f(w)} \right| \leq \frac{\|f\|_{X_{\lambda}}}{\lambda^{\ell}}  \left( \left\|  (\emph{Id}_{n \times n} + L)^{\ell}   \sum_{w \in W}{a_w \delta_w}  \right\|^2_{L^2} - \frac{1}{n}\right)^{\frac{1}{2}}.$$
\end{thm}
The proof shows that the inequality is essentially close to sharp if there is a spectral gap close to $\lambda$. If there is no
spectral gap, there is little reason to assume that the inequality is sharp. If one has access to the eigenfunctions of the Laplacian,
then it is possible to use a much stronger inequality discussed below in \S 2.4. We observe that the inequality only holds for
$f \in X_{\lambda}$: this does \textit{not} mean that the integration scheme will perform poorly on functions outside this function
class. Indeed, it is a fairly common technique to design quadrature rules for a certain class of very smooth functions -- they will
usually still be somewhat effective on functions that are not as smooth as the functions used in the design of the rule. This is,
of course, the main idea in the construction of spherical designs \cite{delsarte, sobolev, stein1} but also appears in other contexts:
the Simpson rule for the integration of a real function on an interval is designed to be exact on quadratic polynomials.

\subsection{Geometric interpretation: Heat Ball Packing.}
Altogether, this suggests that we should use the quantity on the right-hand side in the inequality in Theorem 1, depending only on the set $W$, the weights $a_w$ and the free parameter $\ell$ (but not on the function $f$) as a guideline for how to construct the quadrature rule. This motivates studying the minimization problem
$$ \min_{W \subset V \atop |W| = k} \min_{a_w} \left\|  (\mbox{Id}_{n \times n} + L)^{\ell}   \sum_{w \in W}{a_w \delta_w}  \right\|_{L^2}^2 \qquad \quad \mbox{subject to} ~\sum_{w \in W}{a_w} = 1.$$
Note that we do not actually need to find the precise minimizer of the solution, indeed, we just want the quantity to be small. The topological structure of the minimization problem suggests that there
will be many local minimizer that are not global: this has been empirically observed in related problems in the continuous setting (see e.g. \cite{er}) and there is little reason to assume the discrete
case will be different. Moreover, if the continuous analogies continue to hold, then the actual numerical value of minimizing configuration can be expected to be quite close to that of almost-minimizers which is good news: there are many almost-minimizers (thus easy to find) and they are almost as good as the global minimizer.
$(\mbox{Id}_{n \times n} + L)^{\ell} \delta_w$ is the probability distribution of a random walker starting in $w$ after $\ell$ jumps which lends itself to a geometric interpretation.
\begin{quote}
\textbf{Guideline.} If we manage to find a placement of vertices $W$ with the property that the random walks, weighted by $a_w$, overlap very little, then we have found a good quadrature rule on the graph $G=(V,E)$.
\end{quote}
We observe that this can be reinterpreted as a `packing problem for heat
balls'. This principle has already appeared naturally in the continuous setting
of Riemannian manifolds in work of the second author \cite{stein}. It is even meaningful if we can only sample in a single vertex: we should pick the most central
vertex and that is equivalent to heat diffusing quickly to many other points. In the ideal setting with perfectly distributed heat balls all the weights would be identical (as can be seen in examples with lots of symmetry, see \cite{stein1}); this nicely mirrors an unofficial guideline in numerical integration stating that how far the weights deviate from constant weights can be used as measure of quality and stability of the method.
We summarize that
\begin{enumerate}
\item it is desirable to construct $W \subset V$ equipped with weights $a_w$ such that random walks, starting in $w \in W$ and weighted by $a_w$, intersect each other as little as possible.
\item if we are \textit{not} allowed to chose $W$, we can still use the procedure above to find weights that yield a better result than just uniform sampling.
\end{enumerate}

We emphasize that we do not attempt to solve the heat ball packing problem here -- nor do we expect it to be easily solvable at this level of generality. The main contribution of this paper is to introduce the heat ball packing problem as a fundamental issue with implications for sampling. 

\begin{quote}
\textbf{Problem.} How does one find effective heat ball packing configurations quickly? Are there algorithms leading to effective almost-minimizing configurations? What theoretical guarantees can be proven?
\end{quote}

\subsection{Another interpretation.}
There is a nice bit of intuition coming from a more abstract perspective: the reason why the above method works is that the heat propagator $e^{t\Delta}$ is self-adjoint. More precisely, let us assume $(M,g)$
is a compact manifold normalized to $\mbox{vol}(M) = 1$ and let $\mu$ be a probability measure (our quadrature rule) on it. Smooth functions $f$ have the property that they do not change substantially if we apply the heat propagator for a short time $t>0$ (this is one way of quantifying smoothness) and therefore
$$ \left\langle f, \mu \right\rangle \sim \left\langle e^{t\Delta} f, \mu \right\rangle = \left\langle  f, e^{t\Delta} \mu \right\rangle.$$
In order for this quadrature rule to be effective, we want that $e^{t\Delta} \mu$ is very close to the Lebesgue measure $dx$. However, the heat propagator preserves the $L^1-$mass and thus
$$ 1 = \| \mu\|_{} = \| e^{t\Delta} \mu \|_{L^1} \leq  \| e^{t\Delta} \mu \|_{L^2}.$$
The Cauchy-Schwarz inequality is only sharp if $e^{t\Delta} \mu$ coincides with the constant function and thus minimizing $\| e^{t\Delta} \mu \|_{L^2}$ is a natural way to obtain good quadrature rules. If $\mu$ is a weighted sum of Dirac measures, then $e^{t\Delta}$ turns this, roughly, into a sum of Gaussians (`heat balls') and finding locations to minimize the $L^2-$norm becomes, essentially, a packing problem. We observe that the argument does not single out $L^2$ and minimizing $\| e^{t\Delta} \mu\|_{L^p}$ for $p>1$ will lead to a very similar phenomenon -- it remains to be seen whether there is any advantage to that perspective since $L^2$ is usually easiest to deal with in practice.

\subsection{The Laplacian: a second method.} We believe that our notion of Laplacian may be quite useful in practice (it already appeared in \cite{linderman}). It combines the two desirable properties of
\begin{enumerate}
\item having a symmetric matrix (and thus orthogonal eigenvectors)
\item inducing a diffusion operator that preserves the mean value of the function.
\end{enumerate}
The Kirchhoff matrix $L_1 = D-A$ has the second property but not the first; the normalized Laplacian $L_2 = \mbox{Id}_{n \times n} -  D^{-1/2}AD^{-1/2}$ satisfies the first property but not the second. We believe that, for this reason alone, our notion of a Laplacian may be useful in other contexts as well. We observe that if we can compute eigenvectors directly, then there is a very direct way of approaching the problem directly.

\begin{proposition} Let $W \subset V$ be equipped with weights $a_w$ summing to 1.  Then, for all $0 < \lambda < 1$,
$$ \sup_{f \in X_{\lambda} \atop f \neq 0}  \frac{1}{\|f\|_{X_{\lambda}}} \left| \frac{1}{n} \sum_{v \in V}{f(v)} - \sum_{w \in W}{a_w f(w)} \right| =  \left\|   \frac{1}{n} - \sum_{w \in W}{a_w \delta_w} \right\|_{X_{\lambda}}$$
\end{proposition}

This statement follows easily from $L^2-$duality and will be obtained in the proof of the Theorem as a by-product. The result is true in general but, of course, one cannot compute the quantity on the right-hand side unless one has access to the eigenvectors of the Laplacian. If one has indeed access to either the Laplacian eigenfunctions or at least some of the leading eigenvectors $\phi_1, \dots, \phi_k$ where $k \ll n$ is chosen so that $\lambda_k \sim \lambda$ (simpler to obtain in practice), then optimizing the functional in $W \subset V$ and weights $a_w$ is essentially equivalent to finding good quadrature points. This simple observation is very effective in practice, we refer to numerical examples below.

\section{How to use the Theorem}
The result discussed above has the nice effect of cleanly separating the problem of numerical integration on a graph from the actual graph structure: the geometry of the graph is encoded implicitly 
in the geometry of the random walk. This has the nice effect of providing a uniform treatment but, as a downside, does not provide an immediate method on how to proceed in particular instances. The purpose of this
section is to comment on various aspects of the problem and discuss approaches. Recall that our main result can be written as
$$ \forall~f \in X_{\lambda} \quad   \left| \frac{1}{n} \sum_{v \in V}{f(v)} - \sum_{w \in W}{a_wf(w)} \right| \leq  \|f\|_{X_{\lambda}} \min_{\ell \in \mathbb{N}} \frac{1}{\lambda^{\ell}}  \left( \left\|  (\mbox{Id}_{n \times n} + L)^{\ell}   \sum_{w \in W}{a_w \delta_w}  \right\|_{L^2}^2 - \frac{1}{n}\right)^{\frac{1}{2}}$$

\subsection{The parameters $\lambda$ and $\ell$.} A priori we have no knowledge about the degree of smoothness of the function
and have no control over it. The parameter $\ell$, on the other hand, is quite important and has a nontrivial impact on the minimizing energy configurations for the quantity on the right-hand side. In practice, we have to choose $\ell$ without knowing $\lambda$ (fixing an allowed number of sample points implicitly fixes a scale for $\lambda$). We propose the following basic heuristic.

\begin{quote}
\textbf{Heuristic.} If $\ell$ is too small, then there is not enough diffusion and heat balls interact strongly with themselves. If $\ell$ is too large, the exponentially increasing weight $\lambda^{-\ell}$ is too large. There is an intermediate regime when heat balls starting to interact with nearby heat balls.
\end{quote}
The heuristic is accurate if the graph is locally (at the scale of typical distance between elements of $W$) close to Euclidean. 
On general graphs, the optimal scale of $\ell$ might be more nontrivial to estimate -- we have observed that, in practice,
there is a wide range of $\ell$ yielding good results.

\subsection{The placement of $W$.} Naturally, to avoid intersections of random walkers, we want to place the elements of $W$ as far from each other
as possible. In particular, if the graph is close to having a Euclidean structure, we would expect fairly equi-spaced points to do well. A method that was used in 
\cite{stein1} is to start with a random set of $k$ vertices $\left\{v_1, v_2, \dots, v_k\right\}$ and compute the 
$$ \mbox{total mutual distance} = \sum_{i,j =1}^{k}{d(v_i, v_j)},$$
where $d$ is a metric on $G=(V,E)$. The algorithm then goes through all the vertices and checks whether moving one of them to a neighboring vertices
increases the total mutual distance and, if so, moves the vertex. This is repeated as long as possible.  
The simple numerical examples in \cite{stein1} all have edges with equal weight and the standard combinatorial graph distance
can be used; the optimal type of metric could strongly depend on the type of graph under consideration. 
There is a fairly natural reason why an algorithm of this type has the ability to successfully produce sets of vertices that are very well spread out. We quickly return to the sphere $\mathbb{S}^d$ where a particularly
spectacular justification exists. Let $\sigma$ be the normalized measure on $\mathbb{S}^d$. Then, for any set $X=\left\{x_1, \dots, x_n\right\} \subset \mathbb{S}^d$, we have Stolarsky's invariance principle \cite{dai, stol}
$$   \int_{\mathbb{S}^d} \int_{\mathbb{S}^d} {\|x-y\| d\sigma(x) d\sigma(y)} - \frac{1}{n^2}\sum_{i,j=1}^{n}{\|x_i - x_j\|} = c_d \left[ D_{L^2,\mbox{\tiny cap}}(X) \right]^2,$$
where the quantity on the right is the $L^2-$based spherical cap discrepancy and $c_d$ is a constant only depending on the dimension. The $L^2-$based spherical cap discrepancy is a measure that has been studied in its own right: if the points are evenly distributed, then it is small. This may be a somewhat peculiar case. However, it is not difficult to see that on fairly generic manifolds functionals along the lines of
$$ \sum_{i, j=1}^{n}{\frac{1}{d(x_i, x_j)^{\alpha}}}  \rightarrow \mbox{min} \qquad \mbox{or} \qquad \sum_{i, j=1}^{n}{e^{-\alpha d(x_i, x_j)}}  \rightarrow \mbox{min}$$
converge to the uniform distributions if the number of points is large. Moreover, and this is particularly useful, these types of functionals tend to produce minimal energy configurations
that only weakly depend on the functional being used. On two-dimensional manifolds, the hexagonal lattice seems to be particularly universal (see \cite{blanc}). 
We do not know what kind of interaction functional is optimal on graphs. In practice, one would like to have fast and reliable algorithms that scale well and this seems like a problem of substantial interest.

\subsection{The weights $a_w$.} Once we are given a set $W \subset V$ and a parameter $\ell$, the optimization of the weights is completely straightforward. Observe that
$$ \left\|  (\mbox{Id}_{n \times n} + L)^{\ell}   \sum_{w \in W}{a_w \delta_w}  \right\|_{L^2}^2 = \sum_{w_1,w_2 \in W}{ a_{w_1} a_{w_2} \left\langle   (\mbox{Id}_{n \times n} + L)^{\ell} \delta_{w_1},   (\mbox{Id}_{n \times n} + L)^{\ell} \delta_{w_2}\right\rangle }.$$
This is merely a quadratic form indexed by a $|W| \times |W|$ matrix -- we thus need to solve the semidefinite program
$$  \sum_{w_1,w_2 \in W}{ a_{w_1} a_{w_2} \left\langle   (\mbox{Id}_{n \times n} + L)^{\ell} \delta_{w_1},   (\mbox{Id}_{n \times n} + L)^{\ell} \delta_{w_2}\right\rangle } \rightarrow \mbox{min} \qquad \mbox{subject to} \quad \sum_{w \in W}{a_w} = 1.$$ 
These weights $a_w$ play an important role in further fine-tuning an existing set of vertices $W \subset V$. This is the reason why minimizing the functional can be used to find appropriate weights for \textit{any} given set of vertices $W$: if two vertices in $W$ happen to be somewhat close to each
other, then the quadratic form will take this into account when distributing the weights. Conversely, if one the vertices in $W$ is surprisingly isolated from the other vertices in $W$, then the quadratic form will
increase the weight assigned to that vertex. This is exactly how things should be: points that are oversampling a region in the graph should be given a smaller weight whereas
isolated vertices cover a wider range and are, correspondingly, more important. We refer to \S 4 for numerical examples illustrating this point.

\subsection{Related results.} Sampling on graphs is a fundamental problem and a
variety of approaches have been discussed in the literature \cite{survey}. Sampling is usually done for
the purpose of compression or visualization and not numerical
integration (in particular, vertices are usually not equipped with weights). Our
approach seems to be very different from anything that has been proposed.
The closest construction seems to be a line of research using biased random
walks \cite{gjoka, jin, rasti} based on the idea of
sending random walkers and accepting the vertices they traverse with certain
biased probabilities as sampling points. In contrast, we select points
so that random walkers starting thereavoid each other. Other results seem related in spirit \cite{bermanis}. Our approach is motivated by a
recent approach \cite{stein} to the study of spherical $t-$designs of
Delsarte, Goethals \& Seidel \cite{delsarte} (and Sobolev \cite{sobolev}).
These are sets of points defined on $\mathbb{S}^d$ with the property that they
integrate a large number of low-degree polynomials exactly (we refer to a survey of Brauchart \& Grabner).
\cite{brau}. The second author recently extended some of these results to
weighted points and general manifolds \cite{stein}. These ideas are shown to
be an effective source of good quadrature points in the Euclidean setting in a
paper of Lu, Sachs and the second author \cite{lu}.

\begin{figure}[h!] \label{fig:one}
\begin{minipage}[r]{.4\textwidth}
\begin{tikzpicture}[scale=0.7]
  \tikzstyle{every node}=[circle,inner sep=0pt,minimum size=0.5cm]
    \foreach \y[count=\a] in {10,9,4}
      {\pgfmathtruncatemacro{\kn}{120*\a-90}
       \node at (\kn:3) (b\a) {\small \y};}
    \foreach \y[count=\a] in {8,7,2}
      {\pgfmathtruncatemacro{\kn}{120*\a-90}
       \node at (\kn:2.2) (d\a) {\small \y};}
    \foreach \y[count=\a] in {1,5,6}
      {\pgfmathtruncatemacro{\jn}{120*\a-30}
       \node at (\jn:1.5) (a\a) {\small \y};}
    \foreach \y[count=\a] in {3,11,12}
      {\pgfmathtruncatemacro{\jn}{120*\a-30}
       \node at (\jn:3) (c\a) {\small \y};}
  \draw[dashed] (a1)--(a2)--(a3)--(a1);
  \draw[ultra thick] (d1)--(d2)--(d3)--(d1);
  \foreach \a in {1,2,3}
   {\draw[dashed] (a\a)--(c\a);
   \draw[ultra thick] (d\a)--(b\a);}
   \draw[ultra thick] (c1)--(b1)--(c3)--(b3)--(c2)--(b2)--(c1);
   \draw[ultra thick] (c1)--(d1)--(c3)--(d3)--(c2)--(d2)--(c1);
   \draw[dashed] (b1)--(a1)--(b2)--(a2)--(b3)--(a3)--(b1);
\end{tikzpicture}
\end{minipage} 
\begin{minipage}[l]{.4\textwidth}
\begin{center}
  \begin{tikzpicture}[scale=0.7]
\foreach \a in {1,2,...,24}{
\filldraw (\a*360/24: 3cm) circle (0.09cm);
};
\foreach \a in {1,2,...,24}{
\draw [thick] (\a*360/24: 3cm) --  (\a*360/24 + 360/24: 3cm);
};
\draw [thick] (1*360/24: 3cm) -- (8*360/24: 3cm);
\draw [thick] (2*360/24: 3cm) -- (19*360/24: 3cm);
\draw [thick] (3*360/24: 3cm) -- (15*360/24: 3cm);
\draw [thick] (4*360/24: 3cm) -- (11*360/24: 3cm);
\draw [thick] (5*360/24: 3cm) -- (22*360/24: 3cm);
\draw [thick] (6*360/24: 3cm) -- (18*360/24: 3cm);
\draw [thick] (7*360/24: 3cm) -- (14*360/24: 3cm);
\draw [thick] (9*360/24: 3cm) -- (21*360/24: 3cm);
\draw [thick] (10*360/24: 3cm) -- (17*360/24: 3cm);
\draw [thick] (13*360/24: 3cm) -- (20*360/24: 3cm);
\draw [thick] (16*360/24: 3cm) -- (23*360/24: 3cm);
\draw [thick] (4*360/24: 3cm) circle (0.4cm);
\draw [thick] (7*360/24: 3cm) circle (0.4cm);
\draw [thick] (8*360/24: 3cm) circle (0.4cm);
\draw[thick] (11*360/24: 3cm) circle (0.4cm);
\draw [thick] (16*360/24: 3cm) circle (0.4cm);
\draw[thick] (19*360/24: 3cm) circle (0.4cm);
\draw[thick] (20*360/24: 3cm) circle (0.4cm);
\draw [thick] (23*360/24: 3cm) circle (0.4cm);
   \end{tikzpicture}
\end{center}
\end{minipage} 
\caption{(Left:) The Icosahedron integrates all polynomials up to degree 5 on $\mathbb{S}^2$ exactly (this space is 36-dimensional). (Right:) a subset of 8 vertices integrates 21 of 24 eigenvectors of the McGee graph exactly. Such examples require extraordinary amounts of symmetry and are not generic.}
\end{figure}
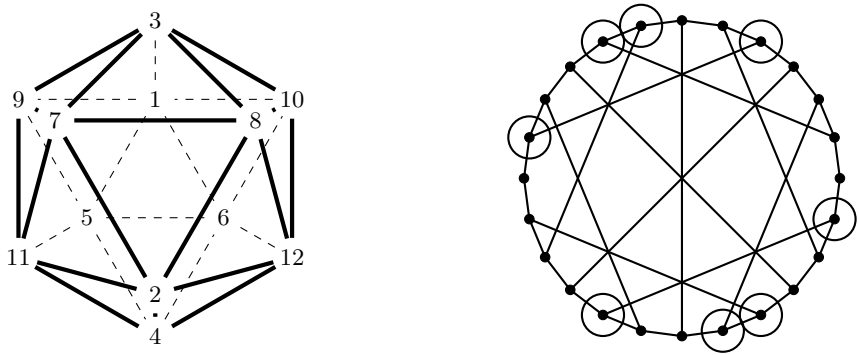

The second author recently
proved \cite{stein1} a generalized Delsarte-Goethals-Seidel bound for graph designs (the analogue of spherical $t-$designs on combinatorial graphs). The main condition in that paper is algebraic (exact integration of
a certain number of Laplacian eigenvectors) as opposed to quantitative (small integration error). \cite{stein1} shows a number of truly remarkable 
quadrature rules on highly structured graphs that were found by numerical search: one of these rules is depicted in Figure 1 and manages with only 8 evaluations to integrate 21 out of a total of 24 eigenvectors exactly. However, these examples are very non-generic, the consequence
of strong underlying symmetries and not likely to be representative of what can be achieved in a typical setting.

\section{Numerical Examples}

\subsection{Importance of weights.} We start with a toy example shown in the Figure below to illustrate the importance of the weights to counterbalance bad geometric distributions. We assume we are given two clusters and an uneven distribution of six sampling points: five end up in one cluster while the other cluster contains a total of five points. We construct a graph based on nearest neighbor distances weighted with a Gaussian kernel. The spreading of heat is rather uneven if all points are given equal weight. Conversely, by adjusting the weight to minimize the $L^2-$norm (constraint to summing to 1 and being nonnegative), a much more balanced distribution is achieved.
\begin{figure}[h!]
\includegraphics[width=0.9\textwidth]{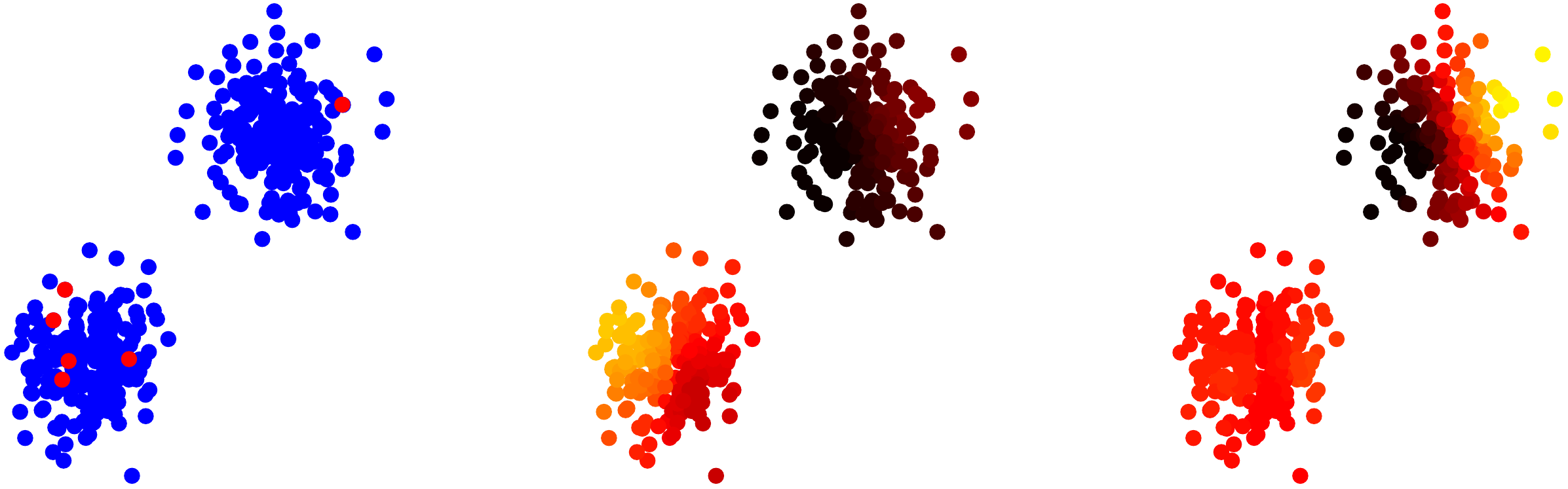}
\caption{Six points selected unevenly in two clusters (left), the heat flow emanating from weighting them all evenly (middle) and the optimal coefficients (0.11, 0.24, 0.05, 0.00, 0.18 and 0.43) for the heat ball packing problem (right).}
\end{figure}

The weights show that, in particular, one point is given weight 0 and another point is given a rather small weight (0.05). This is to counterbalance the clustering of points. The isolated point is given almost half the weight (0.43). We see that the heat distribution in the second cluster is still highly uneven: this shows that it would be preferable to pick another point since the single sampling point is actually quite far from the center of the cluster: if it was closer to the center, it would have received an even larger weight. 

\begin{figure}[h!] \label{fig:one}
\begin{minipage}[r]{.4\textwidth}

\begin{center}
\begin{tikzpicture}[scale=1.1]
\filldraw (0,0) circle (0.05cm);
\draw [thick] (0,0) -- (1,1);
\filldraw (1,1) circle (0.05cm);
\draw [thick] (1,1) -- (1,2);
\draw [thick] (1,1) -- (2,2);
\filldraw (1,2) circle (0.05cm);
\draw [thick] (2,2) -- (3,2) -- (3,3) -- (2,3) -- (2,2);
\filldraw (2,2) circle (0.05cm);
\filldraw (3,2) circle (0.05cm);
\draw [thick] (3,3) -- (3,4) -- (2,4) -- (1,4) -- (2,2);
\filldraw (3,3) circle (0.05cm);
\filldraw (2,3) circle (0.05cm);
\filldraw (1,4) circle (0.05cm);
\filldraw (2,4) circle (0.05cm);
\filldraw (3,4) circle (0.05cm);
\draw [thick] (1,2) -- (2,3) -- (2,4) -- (3,2);
\draw [thick] (2,3) -- (1,4);
\draw [thick] (1,1) circle (0.3cm);
\node at (0.1, 1) {0.19};
\draw [thick] (3,3) circle (0.3cm);
\node at (4, 3) {0.81};
\draw [thick] (2,4) circle (0.3cm);
\node at (2, 4.7) {0};
\end{tikzpicture}
\end{center}

\end{minipage} 
\begin{minipage}[l]{.4\textwidth}

\begin{center}
\begin{tikzpicture}[scale=1.1]
\filldraw (0,0) circle (0.05cm);
\draw [thick] (0,0) -- (1,1);
\filldraw (1,1) circle (0.05cm);
\draw [thick] (1,1) -- (1,2);
\draw [thick] (1,1) -- (2,2);
\filldraw (1,2) circle (0.05cm);
\draw [thick] (2,2) -- (3,2) -- (3,3) -- (2,3) -- (2,2);
\filldraw (2,2) circle (0.05cm);
\filldraw (3,2) circle (0.05cm);
\draw [thick] (3,3) -- (3,4) -- (2,4) -- (1,4) -- (2,2);
\filldraw (3,3) circle (0.05cm);
\filldraw (2,3) circle (0.05cm);
\filldraw (1,4) circle (0.05cm);
\filldraw (2,4) circle (0.05cm);
\filldraw (3,4) circle (0.05cm);
\draw [thick] (1,2) -- (2,3) -- (2,4) -- (3,2);
\draw [thick] (2,3) -- (1,4);
\draw [thick] (1,1) circle (0.3cm);
\node at (0.1, 1) {0.2};
\draw [thick] (3,3) circle (0.3cm);
\node at (4, 3) {0.4};
\draw [thick] (1,4) circle (0.3cm);
\node at (1, 4.7) {0.4};
\end{tikzpicture}
\end{center}
\end{minipage} 
\caption{Two optimal configurations for $\ell=3$ on 3 vertices.}
\label{twooptimal}
\end{figure}
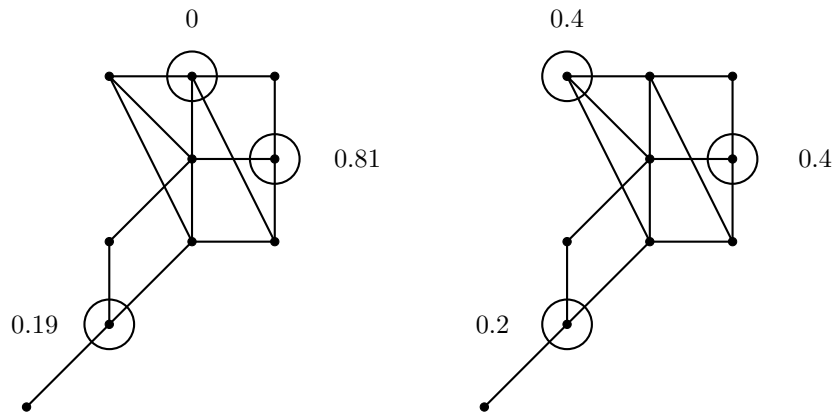

Figure \ref{twooptimal} shows an example on a small graph with 10 vertices and
15 edges. More precisely, we optimize $\|(\mbox{Id}_{10 \times 10} + L)^3
\sum_{w} a_w \delta_w\|_{L^2}$ over all sets $W$ with three vertices. In the first
example, we see that there is one very central node that distributes well
throughout the network, another weight is actually set to 0. This is a
consequence of constraining the optimization to non-negative weights $a_w \geq 0$. This
constraint is not required by the Theorem, but makes the optimization easier,
and is well-motivated by classical methods in numerical integration.  If we move
the point that was assigned weight 0, then the weight splits evenly (the value
of the functional barely changes).

\subsection{MNIST}
Our explicit example is as follows: we consider the data set MNIST, a
collection of handwritten digits represented as $28 \times 28$ pixels (each
pixel being either 0 or 1). For simplicity, we only consider the subset
comprised of the digits 0 and 1. The challenge problem is to figure out the
proportion of elements in the set that are handwritten digits that are 1's
(which is one half). This ties in to our example in the beginning: suppose we
did not know the precise proportion of 1's and the data is unlabeled. Labeling
the data is expensive: the function evaluation would be one human being looking
at a picture and labeling it, which is costly.
\begin{figure}[h!]
\begin{tikzpicture}
\node at (0,0){\includegraphics[width=0.7\textwidth]{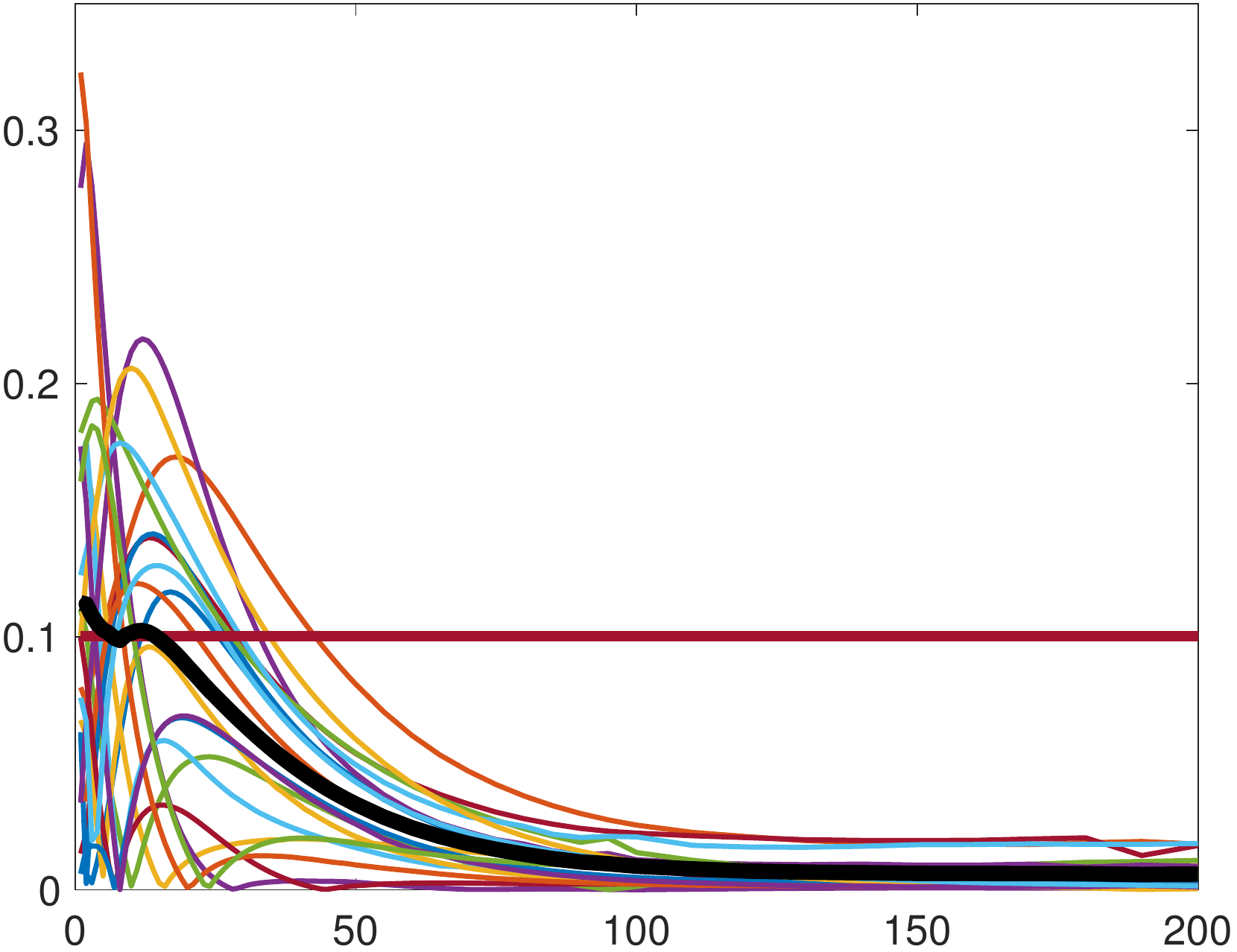}};
\node at (0,-4.5) {Number $\ell$ of diffusion steps.};
\node [rotate=90] at (-5.6,0) {Relative error};
\end{tikzpicture}
\label{fig:eigenvecs_bound}
\caption{Estimating digits in MNIST: the relative integration error for 20
	different sets of points of size 50 and how it evolves as a function
	depending on $\ell$ (leading to different selection of weights). The
	average of these 20 curves is shown as the bold black line, sampling error
	for randomly chosen points is the red line.}
\end{figure}
However, these pictures are merely $\left\{0,1\right\}-$vectors in
$\mathbb{R}^{784}$. As is commonly done, we reduce the dimensionality of the
data by projecting onto its first ten principal components. We build a graph by
connecting every element to its 10-nearest neighbors (in Euclidean distance)
weighted with a Gaussian kernel and then symmetrize by averaging the
resulting adjacency matrix with its transpose.  It is reasonable to assume that the
indicator function of 1's is smooth over a graph defined by that notion of
distance: handwritten digits looking like a 1 should be close to other
handwritten digits that look like 1. 
We then proceed as outlined above: we sample random points, move them iteratively so that they are far away from each other and then adjust weights by solving the semidefinite program. The result is plotted against the parameter $\ell$ and compared to uniform weights on random points (red); the picture shows 20 different sets of points, the evolution of their integration error depending on $\ell$ as well as their average (black). We observe that
 for the right parameter range of $\ell$, the obtained numerical integration scheme performs much better but the precise performance depends on the points chosen. This highlights the need for fast, stable and guaranteed ways of approximately solving the heat ball packing problem.

\subsection{Using Eigenvectors of Laplacian} 
This section studies the same example as above, estimating the proportion of handwritten digits `1' in MNIST, but assumes additionally that we are able to access the eigenvectors of the Laplacian associated to the largest few eigenvalues exactly. We set $\lambda = 0.994$ close to 1 leading to a space $X_{\lambda}$ spanned by very few of the smoothest eigenvectors, sample random points, move them far apart and make use of 
$$ \sup_{f \in X_{\lambda} \atop f \neq 0}  \frac{1}{\|f\|_{X_{\lambda}}} \left| \frac{1}{n} \sum_{v \in V}{f(v)} - \sum_{w \in W}{a_w f(w)} \right| =  \left\|  -\frac{1}{n} + \sum_{w \in W}{a_w \delta_w} \right\|_{X_{\lambda}}$$
to explicitly optimize the weights. 
\begin{figure}[h!]
\begin{tikzpicture}
\node at (0,0) {\includegraphics[width=0.6\textwidth]{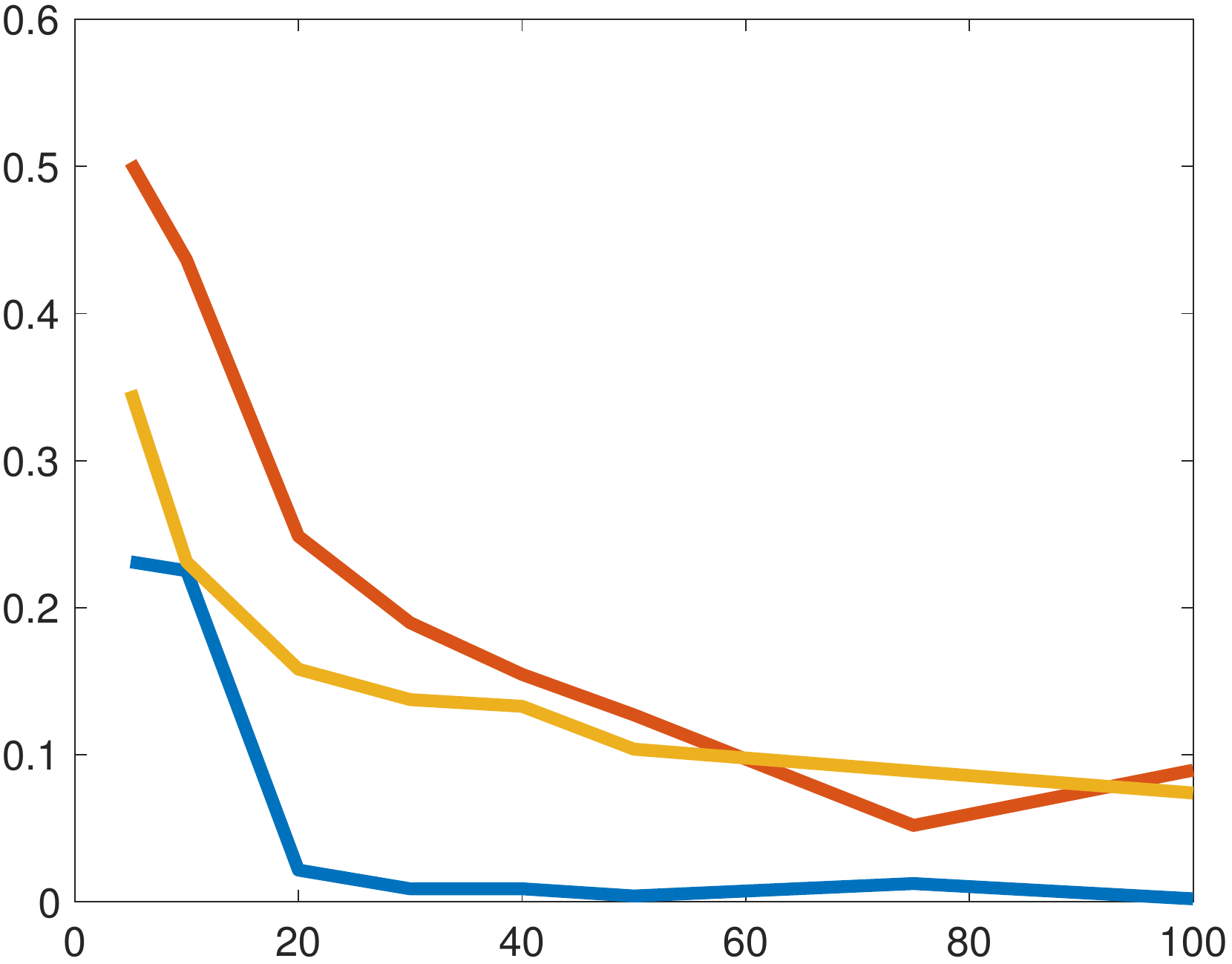}};
\node [rotate=90] at (-5,0) {Relative error};
\node at (0,-4) {Sampling size};
\end{tikzpicture}
	\caption{Numerical Integration with access to eigenvectors: our method (blue) compared to sampling in random points (yellow) and sampling in the location in which we move the points (red).}
\label{fig:eigenvecs_performance}
\end{figure}
Our method (blue) is shown to perform exceedingly well; the weights are crucial, sampling over the same points (red) is even worse than sampling over random points (yellow) which decays as $(\#\mbox{size of subset})^{-1/2}$  (averaged
	over 200 random samples).
As mentioned above, having direct access to the eigenvectors implies that the bound is often close to sharp. This is illustrated in the subsequent Figure 6 where we replace the indicator function of the `1'-digits in the subset of MNIST comprised of digits 0 and 1 by its mollification obtained from projecting onto the first 6 eigenvectors. We optimize in the weights and obtain
$$  \left\|\frac{1}{n} -   \sum_{w \in W}{a_w \delta_w} \right\|_{X_{\lambda}} \qquad \mbox{as an upper bound (red)}$$
that is then compared to the error on the function $f$ (blue). We observe that as soon as the sample size exceeds a certain limit, integration becomes exact. It is quite desirable to obtain a better understanding of the interplay of parameters involved: suppose we are given a set of $k$ well-distributed points and optimize the weights so as to minimize $\left\|  -n^{-1} + \sum_{w \in W}{a_w \delta_w} \right\|_{X_{\lambda}}$, what is the interplay between $k$, the $X_{\lambda}$ space and the performance of the arising quadrature rule? Or, put differently, how does the integration error in a $X_{\mu}$ space depend on the space $X_{\lambda}$ that was used to determine the weights? This question is clearly of great relevance in applications.
\begin{figure}[h!]
\begin{tikzpicture}
\node at (0,0){\includegraphics[width=0.6\textwidth]{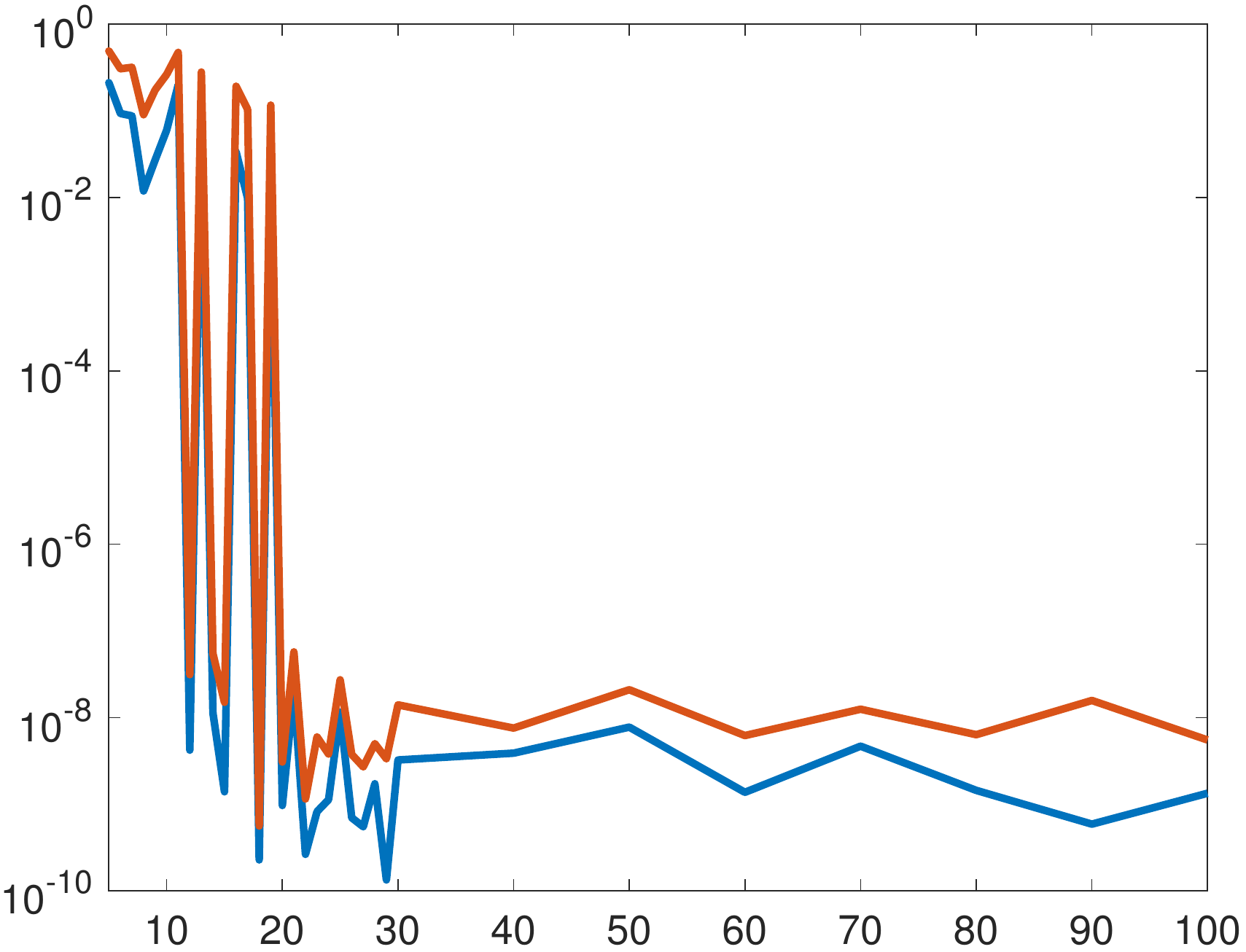}};
\node [rotate=90] at (-5,0) {Integration error};
\node at (0,-4) {Sampling size};
\end{tikzpicture}
	\caption{
Error for the smoothed indicator function on the 1's (blue) and the theoretical upper bound on the integration error (red). It depends strongly on $W$ and, after
some initial fluctuation settles down to essentially exact integration; the theoretical upper bound matches the performance on the particular instance.}
	%\lambda = 0.994
\label{fig:eigenvecs_bound}
\end{figure}

The phenomenon, which seems to be generic and easily observed in most examples, is illustrated in Figure 7.  We downsampled the MNIST dataset of digits `0' and `1' to consist of a total of 1000 points, and constructed the graph as before. We then choose a subset of points $W$ of size $100$, increase their mutual distance and fix them as quadrature points. 
\begin{figure}[h!]
\begin{tikzpicture}
\node at (0,0){\includegraphics[width=0.6\textwidth]{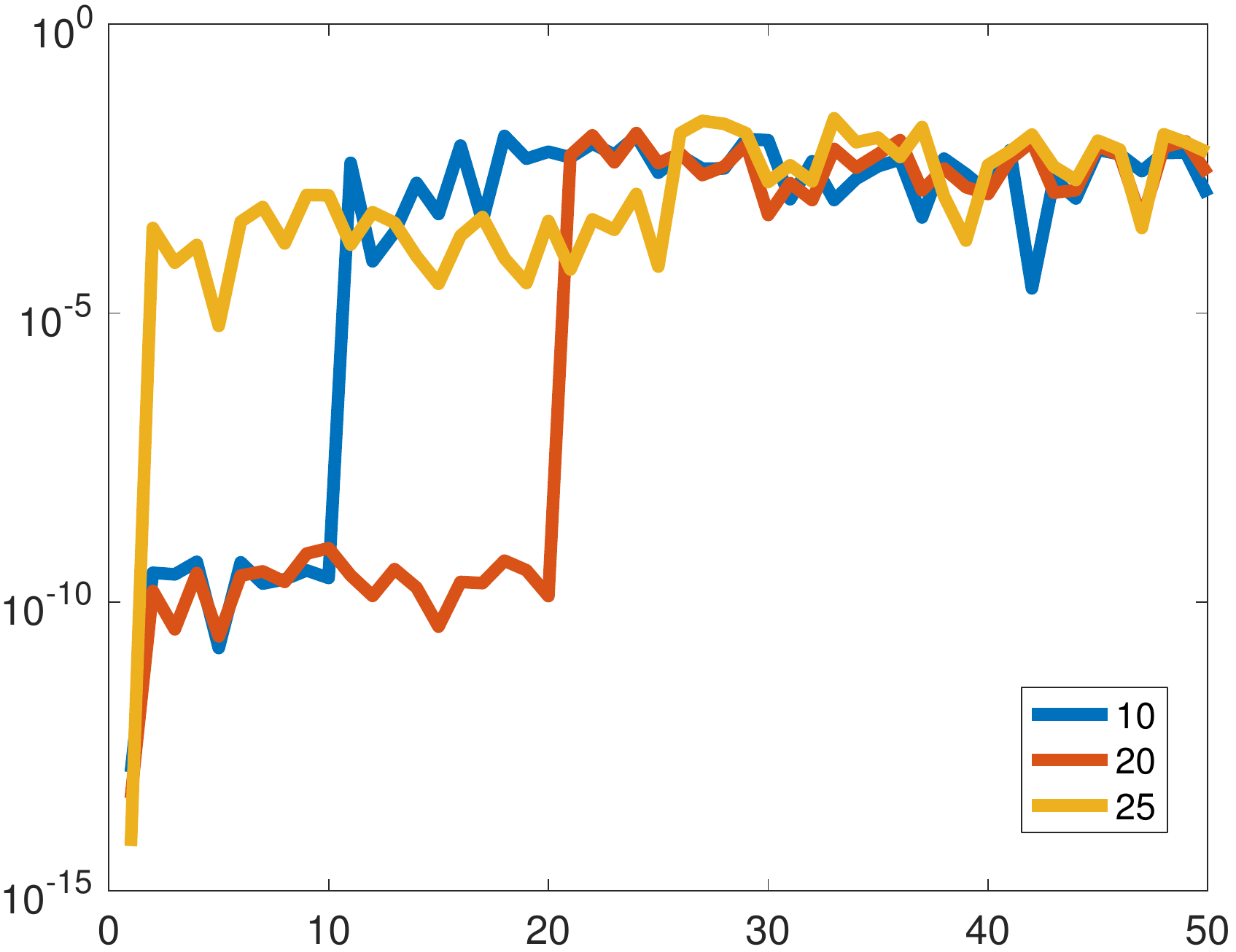}};
\node [rotate=90] at (-5,0) {Integration error};
\node at (0,-4) {Eigenvector};
\end{tikzpicture}
\vspace{-5pt}
	\caption{
	Integration error of three quadrature rules with weights fine-tuned in three different $X_{\lambda}-$spaces on the  first 50 eigenfunctions. 
}
	%\lambda = 0.994
\label{fig:eigenvecs_bound}
\end{figure}

Finally, we optimize their weights in three different $X_{\lambda}$ spaces where $\lambda$ is chosen such that the dimensions of the spaces are 10, 20 and 25 (i.e. they contain the first 10, 20 and 25 eigenfunctions, respectively). We then plot the integration error of these three quadrature rules on the first 50 eigenfunctions. If we optimize the weights according to the $X_{\lambda}$ space containing the first 10 eigenfunctions, then the first 10 eigenfunctions are essentially integrated exactly, the subsequent integration error is small. The same is true for optimization in the space containing the first 20 eigenfunctions. Then the behavior changes abruptly: if we optimize over the first 25 eigenfunctions, then the error on those 25 eigenfunctions is small ($\sim 10^{-4}$) and, as in the other examples, increases afterwards.
This seems to be typical: for any given set $W \subset V$, there seems to be a range of $X_{\lambda}-$spaces such that optimizing parameters leads to exact integration in $X_{\lambda}$. Once their dimension exceeds a certain (sharp) threshold, the error is still small in $X_{\lambda}$ but many orders of magnitude larger than before. This sharp phase transition could serve as another measure of quality of $W$ that may be useful in judging the quality of algorithms finding $W$ (the largest number of eigenfunctions that can be integrated exactly using $W$ for some choice of weights, a measure already studied in \cite{stein1}).

\section{Proof of the Theorem}
\begin{proof}

We write the integration error as the inner product of two vectors and decompose it as
\begin{align*}
  \frac{1}{n} \sum_{v \in V}{f(v)} - \sum_{w \in W}{a_w f(w)}   &=   \sum_{v \in V}{ \left(\frac{1}{n} - \sum_{w \in W}{a_w \delta_w}(v) \right) f(v)}       \\
&= \sum_{k=1}^{n}{ \left\langle  \frac{1}{n} - \sum_{w \in W}{a_w \delta_w} , \phi_k \right\rangle  \left\langle f , \phi_k \right\rangle}.
\end{align*} 
Since $f \in X_{\lambda}$, we have that $\left\langle f, \phi_k \right\rangle = 0$ unless $|\lambda_k|\geq \lambda$. A simple application
of the Cauchy-Schwarz inequality then shows that
	\begin{equation*}\label{eqn:eigenvecs_bound}
		\left| \frac{1}{n} \sum_{v \in V}{f(v)} - \sum_{w \in W}{a_w f(w)} \right| \leq  \left\|  \frac{1}{n} - \sum_{w \in W}{a_w \delta_w} \right\|_{X_{\lambda}} \|f\|_{X_{\lambda}}.
	\end{equation*}
More precisely, $L^2-$duality implies that this step is not lossy since 
$$ \sup_{f \in X_{\lambda} \atop f \neq 0}  \frac{1}{\|f\|_{X_{\lambda}}} \left| \frac{1}{n} \sum_{v \in V}{f(v)} - \sum_{w \in W}{a_w f(w)} \right| =  \left\|  \frac{1}{n} - \sum_{w \in W}{a_w \delta_w} \right\|_{X_{\lambda}}.$$

For any function $g: V \rightarrow \mathbb{R}$, we have
\begin{align*}
 \left\| g \right\|_{X_{\lambda}}^2 = \sum_{|\lambda_k | \geq \lambda}{ \left |\left\langle g, \phi_k \right\rangle\right|^2}  &\leq \frac{1}{\lambda^{2\ell}}  \sum_{k=1}^{n}{|\lambda_k |^{2\ell} \left|\left\langle g, \phi_k \right\rangle\right|^2} \\
&= \frac{1}{\lambda^{2\ell}} \left\| \sum_{k=1}^{n}{ \lambda_k^{\ell} \left\langle g, \phi_k \right\rangle} \phi_k\right\|^2_{L^2}  \\
&= \frac{1}{\lambda^{2\ell}} \left\| (\mbox{Id}_{n \times n} + L)^{\ell} g  \right\|^2_{L^2}.
\end{align*}
We observe that this inequality is also valid if $g \notin X_{\lambda}$ since ever step is a valid bound from above and $\|\cdot\|_{X_{\lambda}}$ is defined on all functions as a semi-norm (we note, however, that if $g \notin X_{\lambda}$, then the inequality will usually be far from sharp). We use this inequality for
$$ g = \frac{1}{n} - \sum_{w \in W}{a_w \delta_w}  $$
to conclude that
$$  \left\|  \frac{1}{n} - \sum_{w \in W}{a_w \delta_w} \right\|_{X_{\lambda}} \|f\|_{X_{\lambda}} \leq  \frac{\|f\|_{X_{\lambda}}}{\lambda^{\ell}} \left\| (\mbox{Id}_{n \times n} + L)^{\ell}  \left(  \frac{1}{n} - \sum_{w \in W}{a_w \delta_w}  \right)  \right\|_{L^2}.$$
We observe that $\mbox{Id}_{n \times n} + L$ is the generator of the diffusion. It is a linear operator for which constant functions are invariants. This implies 
$$  (\mbox{Id}_{n \times n} + L)^{\ell}  \left(    \frac{1}{n} - \sum_{w \in W}{a_w \delta_w} \right) =  \frac{1}{n}  -  (\mbox{Id}_{n \times n} + L)^{\ell}   \sum_{w \in W}{a_w \delta_w}$$
We now show that the operator $(\mbox{Id}_{n \times n} + L)$ preserves the average value of the function. It suffices to show that $L$ maps the average value of a function to 0. We use the definition 
$$ L =  A' - D' = \frac{1}{d_{\mbox{max}}}\left(A - D\right).$$
It thus suffices to show that $A-D$ maps every function to a function with mean value 0. This follows from changing the order of summation
$$ \sum_{i \in V}{( A f)(i)} = \sum_{i \in V}{ \sum_{j \in V}{   e_{ij}    f(j)}} = \sum_{j \in V}{ \sum_{i \in V}{    e_{ij}  f(j)}}  = \sum_{i \in V}{D_{ii} f(i)} = \sum_{i \in V}{(Df)(i)}.$$
This implies that if we normalize the weights so that constants are being integrated exactly, i.e. 
$$ \sum_{w \in W}{a_w} = 1,$$
then the mean value of 
$$  (\mbox{Id}_{n \times n} + L)^{\ell}   \sum_{w \in W}{a_w \delta_w} \qquad \mbox{is exactly} \quad \frac{1}{n}.$$
Squaring out implies
\begin{align*}
\left\|  \frac{1}{n}  -  (\mbox{Id}_{n \times n} + L)^{\ell}   \sum_{w \in W}{a_w \delta_w} \right\|_{L^2}^2 &= \frac{1}{n}+ \left\|  (\mbox{Id}_{n \times n} + L)^{\ell}   \sum_{w \in W}{a_w \delta_w}  \right\|_{L^2}^2 \\ 
& - 2 \left\langle \frac{1}{n},  (\mbox{Id}_{n \times n} + L)^{\ell}   \sum_{w \in W}{a_w \delta_w} \right\rangle \\
&=  \left\|  (\mbox{Id}_{n \times n} + L)^{\ell}   \sum_{w \in W}{a_w \delta_w}  \right\|_{L^2}^2 - \frac{1}{n}.
\end{align*}
 Altogether, we have shown
$$  \sup_{f \in X_{\lambda}}  \left| \frac{1}{n} \sum_{v \in V}{f(v)} - \sum_{w \in W}{f(w)} \right| \leq \frac{\|f\|_{X_{\lambda}}}{\lambda^{\ell}}  \left( \left\|  (\mbox{Id}_{n \times n} + L)^{\ell}   \sum_{w \in W}{a_w \delta_w}  \right\|_{L^2}^2 - \frac{1}{n}\right)^{\frac{1}{2}}.$$
\end{proof}

\textbf{Acknowledgement} GCL was supported by NIH grant \#1R01HG008383-01A1 (PI: Yuval Kluger) and U.S. NIH MSTP Training Grant T32GM007205.

\end{document}